\documentclass{amsart}
\usepackage[left=25mm, right=25mm, top=1in, bottom=1in]{geometry}
\usepackage{amsmath,bm,amssymb}
\usepackage{scrextend} 
\allowdisplaybreaks
\usepackage[normalem]{ulem}
\usepackage[dvipsnames]{xcolor}
\usepackage{setspace}
\usepackage{amsbsy}
\usepackage{hyperref}
\usepackage{graphicx}
\usepackage[overload]{empheq}

\theoremstyle{plain}
\newtheorem{theorem}{Theorem}

\newtheorem{definition}[theorem]{Definition}

\newtheorem{remark}[theorem]{Remark}

\newcommand{\del}{\partial}
\newcommand{\diver}{\operatorname{div}}

\hypersetup{
    colorlinks=true,
    linkcolor=blue,
    filecolor=blue,      
    urlcolor=blue,
}

\numberwithin{equation}{section}


\begin{document}

\title{Analysis of a nonisothermal Maxwell--Stefan system \\
with degenerate thermal conductivity}

\author[S. Georgiadis]{Stefanos Georgiadis}
        \address[S. Georgiadis]{Division of Science and Mubadala Arabian Center for Climate and Environmental Science\\ New York University
        Abu Dhabi \\ United Arab Emirates}
        \email[]{\href{seg577@nyu.edu}{seg577@nyu.edu}}

\begin{abstract}
  We prove the global-in-time existence of weak solutions for the Maxwell--Stefan--Fourier system with physically motivated degenerate thermal conductivity $\kappa(\theta)=\theta^\alpha$, for $0<\alpha\leq2$. In contrast to existing nonisothermal compressible theories based on nondegenerate conductivities, the entropy inequality no longer gives control of the quantities $\nabla\theta$ and $\nabla \log\theta$. The main new ingredient is a renormalized energy estimate, which yields compactness of the temperature and allows the identification of the degenerate heat flux.
\end{abstract}

\maketitle



\section{Introduction}

The Maxwell--Stefan system is a fundamental model for multicomponent mass transport in fluid mixtures. Unlike Fickian diffusion models, which relate species fluxes directly to concentration gradients, the Maxwell--Stefan formulation is based on a balance between thermodynamic driving forces and frictional interactions among the different species. This framework naturally accounts for the coupling between component fluxes and has therefore become a standard tool in the modeling of transport processes in gases, liquids, porous media, membrane separations, and electrochemical devices.

We consider a nonisothermal Maxwell--Stefan system, also known as Maxwell--Stefan--Fourier, in a bounded domain $\Omega \subset \mathbb{R}^3$, consisting of the mass balance equations
\begin{equation}\label{mass}
    \partial_t \rho_i+\diver J_i=0, \qquad i\in\{1,\dots,n\},   
\end{equation}
coupled to an energy equation
\begin{equation}\label{energy}
    \partial_t(\rho e)+\diver J_e=0.    
\end{equation}
The mass fluxes $J_i$ are determined by the Maxwell--Stefan relations
\begin{align}
    - \sum_{j\not=i} \frac{\rho_j J_i-J_j\rho_i}{D_{ij}} &= d_i, \quad i\in\{1,\dots,n\} \label{linearsystem} \\
    \sum_{i=1}^n J_i &= 0 \label{constraint},
\end{align}
where $D_{ij}$ are positive symmetric interaction coefficients and $d_i$ denotes the thermodynamic driving force acting on the $i$-th species. The energy flux $J_e$ is given by
\begin{equation}\label{heatflux}
    J_e = -\kappa(\theta) \nabla\theta + \sum_{j=1}^n\frac{h_j}{\rho_j} J_j,
\end{equation}
where $\kappa$ is the thermal conductivity and $h_i$ is the partial enthalpy of the $i$-th component.

For an ideal-gas mixture with constant heat capacities, the entropy density takes the form
\[
    H(t) = \int_\Omega \sum_{i=1}^n \left( R_i\rho_i \left(\log \rho_i-1\right) - c_i \rho_i \log\theta \right) \,dx, 
\]
where $R_i>0$ is the engineering gas constant and $c_i>0$ is the heat capacity assumed to be constant for ideal gases. The constitutive quantities are then given by
\[
h_i=(R_i+c_i)\rho_i\theta, \qquad d_i=R_i\nabla(\rho_i\theta), \qquad \rho e=\sum_{i=1}^n c_i\rho_i\theta
\]
(see \cite{GJ24} for more details). Throughout this work we assume, for simplicity, $R_i=c_i=1$ for all $i$ and normalize the total density by
\[
\rho=\sum_{i=1}^n\rho_i=1.
\]
Under these assumptions, the second term in \eqref{heatflux} vanishes due to the constraint \eqref{constraint} and the energy flux reduces to the conductive contribution alone. This simplification is intentional: the purpose of the paper is to isolate the influence of the thermal conductivity and, in particular, its degeneracy at low temperatures.

The mathematical analysis of nonisothermal Maxwell--Stefan systems is by now well established. Global existence results are available in the isothermal setting \cite{JS13} and, more recently, in the nonisothermal setting \cite{HJ21, GJ24}. A common feature of all existing theories is the assumption of a nondegenerate thermal conductivity,
\[
    \kappa(\theta) \gtrsim 1+\theta^\alpha,
\]
with $\alpha\ge1$ and typically with $\alpha=2$. This assumption enters the analysis in two related but distinct ways. First, the entropy production contains the heat-dissipation term
\[
\kappa(\theta)\frac{|\nabla\theta|^2}{\theta^2}.
\]
When $\kappa$ is bounded below by a positive constant, this gives the logarithmic temperature estimates
\[
\nabla\theta\in L^2(\Omega\times(0,T)) \quad \textnormal{ and } \qquad\nabla\log\theta\in L^2(\Omega\times(0,T))
\]
which play a fundamental role in the compactness theory.

Second, in order to identify the heat flux in the weak limit, one needs additional temperature estimates obtained from the energy equation. In the nondegenerate theory, these estimates are derived by testing the energy equation with suitable powers of the temperature, for instance with $\theta$ or $\theta^{\alpha-1}$, depending on the growth of $\kappa$. The physically relevant situation is different. In kinetic theory, the thermal conductivity of a hard-sphere gas behaves like
\[
\kappa(\theta)\sim\theta^\frac12,
\]
and therefore degenerates as $\theta\to0$. More generally, we consider the degenerate law
\[
\kappa(\theta)=\theta^\alpha,\qquad 0<\alpha\le2.
\]
Once this physically motivated law is adopted, the entropy production no longer controls $\nabla\log\theta$. It only gives
\[
\nabla\theta^{\frac\alpha2}\in L^2(\Omega\times(0,T))
\]
and consequently, the existing existence theories no longer apply. In fact, this estimate is insufficient to identify the degenerate heat flux
\[
-\theta^\alpha\nabla\theta = -\frac1{\alpha+1}\nabla\theta^{\alpha+1}.
\]
The natural formal estimate would be obtained by testing the temperature equation with $\theta^{\alpha+1}$. However, this test function is not admissible under the available uniform estimates: the regularity needed to use it is precisely the regularity one is trying to prove. Thus the passage from nondegenerate to degenerate heat conductivities is not a simple modification of the existing theory; the standard temperature-estimate step becomes circular.

The main contribution of this paper is to overcome this obstruction by a renormalized and truncated testing argument. Instead of testing directly with $\theta^{\alpha+1}$, we use the bounded renormalizations
\[
\beta_k(\theta)=T_k(\theta^{\alpha+1}),
\]
where $T_k$ denotes a truncation. These renormalizations are admissible at the approximation level, and they yield, after passing to the limit $k\to\infty$, the uniform estimates
\[
\theta\in L^\infty(0,T;L^{\alpha+2}(\Omega)), \qquad \nabla\theta^{\alpha+1}\in L^2(\Omega\times(0,T)).
\]
The resulting identity is obtained in renormalized form and gives an inequality rather than a classical equality, reflecting the limited regularity of the degenerate temperature equation. This estimate provides the compactness of the temperature and permits the identification of the weak limit of the degenerate heat flux.

We emphasize that this estimate is not specific to the algebraic structure of the Maxwell--Stefan fluxes. It is a coercive estimate associated with the degenerate conductive operator
\[
\diver(\theta^\alpha\nabla\theta) = \frac1{\alpha+1}\Delta\theta^{\alpha+1}.
\]
The Maxwell--Stefan structure is essential for the compactness of the species densities and for the passage to the limit in the cross-diffusion system, but the renormalized temperature estimate itself is model-independent. This suggests that the mechanism developed here may also be useful in other heat-conducting systems, in particular in the analysis of compressible Navier--Stokes--Fourier models with degenerate thermal conductivity.

In order to close the compactness argument for the full Maxwell--Stefan--Fourier system, however, we still assume finite initial entropy, in particular $\log\theta^0\in L^1(\Omega)$. This assumption propagates along the approximation and implies that the limiting temperature remains positive almost everywhere. Thus the existence theorem proved in this paper does not produce, nor does it allow within its finite-entropy class, zero-temperature regions of positive measure.

Nevertheless, the degeneracy of the thermal conductivity points to a different mathematical mechanism from the one present in the nondegenerate theory. Indeed, if one isolates the temperature equation, the energy balance reduces to the porous-medium type equation
\[
\partial_t\theta-\frac1{\alpha+1}\Delta\theta^{\alpha+1}=0.
\]
This equation admits Barenblatt self-similar solutions with compact support, at least before interaction with the boundary. Hence the degenerate conductive law itself is compatible with finite-speed propagation and positive-measure regions where the temperature vanishes. In the present theorem such states are excluded by the finite-entropy assumption, not by the degenerate heat equation.

This observation highlights the limitation of the classical ideal-gas entropy used here. With constant heat capacity, the entropy contains the singular contribution $c_v\log\theta$, which diverges as $\theta\to0$. Therefore finite entropy necessarily rules out genuine zero-temperature states. A more realistic low-temperature model would require a temperature-dependent heat capacity satisfying $c_v(\theta)\to0$ as $\theta\to0$, in accordance with the third law of thermodynamics. Such a modification would change the entropy structure and may lead to a weak-solution framework that no longer excludes thermal vacuum. The analysis of this genuinely degenerate entropy regime is left for future work.

To the best of our knowledge, this is the first existence theory for weak solutions to a nonisothermal Maxwell--Stefan system with physically degenerate thermal conductivity \cite{MPZ15, BJPZ22}. In fact, the same phenomenon appears in the theory of compressible Navier--Stokes--Fourier systems. Existing weak-solution theories rely on nondegenerate heat conductivities and logarithmic temperature estimates; see \cite{BD06, AP15, FNS16, FN17} and the references therein. There are indeed some works covering a degenerate thermal conductivity \cite{DGZ17, LPP26} but the analysis is only in one space dimension and for strong solutions. Another work that points to the degenerate direction can be found in \cite{FJSZ21}, however the model under investigation is a density–temperature diffusion system with a substantially different structure from the heat-conducting mixture and fluid models considered in the Maxwell--Stefan and Navier--Stokes--Fourier literature. We are not aware of any multidimensional weak-solution theory for heat-conducting compressible fluids that is compatible with physically degenerate conductivities of the form $\kappa(\theta)=\theta^\alpha$, $\alpha>0$.


\section{Assumptions and main results}

We impose the following assumption:

\begin{itemize}
    \item[(A1)] Domain: $\Omega \subset \mathbb{R}^3$ is a bounded domain with Lipschitz boundary, and the time horizon is $T>0$.
    \item[(A2)] Initial data: $\rho_i(x,0)=\rho_i^0(x) \in L^\infty(\Omega)$ such that $\rho_i^0 \geq 0$ in $\Omega$ for all $i\in\{1,\dots,n\}$ and $\sum_{i=1}^n\rho_i^0=1$; $\theta(x,0) = \theta^0(x) \in L^{\alpha+2}(\Omega)$, such that $\log\theta^0 \in L^1(\Omega)$.
    \item[(A3)] No-flux boundary conditions: $J_i \cdot \nu = 0$ for all $i\in\{1,\dots,n\}$ and $J_e \cdot \nu = 0$ on $\partial\Omega$, where $\nu$ is the normal vector to the boundary. 
    \item[(A4)] Thermal conductivity: $\kappa(\theta)= \theta^\alpha$, for $0<\alpha\leq2$.
    \item[(A5)] Constants: $c_i=R_i=1$ for all $i\in\{1,\dots,n\}$.
\end{itemize}

\begin{definition}[Weak solutions] \label{def}
    A weak solution to \eqref{mass}--\eqref{energy} consists of functions $\rho_i\geq0$, $\theta>0$ almost everywhere, with
    \[ 
    \del_t\rho_i, \del_t\theta \in L^2(0,T;W^{-1,2}(\Omega)), \qquad J_i, J_e \in L^2(\Omega\times(0,T))
    \]
    for all $i \in\{1,\dots,n\}$, where $J_i$ satisfies the linear system \eqref{linearsystem} in $\mathcal{D}'(\Omega\times(0,T))$ and the constraint \eqref{constraint} almost everywhere, and $J_e$ is given by \eqref{heatflux}, such that
    \begin{align}
        & \int_0^T \langle \del_t\rho_i, \phi\rangle \,dt - \int_0^T \int_\Omega J_i \cdot \nabla\phi \, dxdt = 0, \label{weakmass} \\
        & \int_0^T \langle \del_t\theta, \psi\rangle \,dt - \int_0^T \int_\Omega J_e \cdot \nabla\psi \, dxdt = 0, \label{weakenergy}
    \end{align}
    for test functions $\phi,\psi \in L^2(0,T;H^1(\Omega))$, where the brackets $\langle\cdot,\cdot\rangle$ denote the duality pairing ${W^{-1,2}(\Omega)\times H^1(\Omega)}$.
\end{definition}

\begin{theorem}[Existence of weak solutions] \label{thmex}
    Let assumptions (A1)--(A5) hold. There exists a weak solution to \eqref{mass}--\eqref{heatflux} according to Definition \ref{def}, with $\rho_i\geq0$ for all $i\in\{1,\dots,n\}$ and $\theta>0$ almost everywhere, with the additional regularity
    \[
    \sqrt{\rho_i} \in L^\infty(\Omega\times(0,T)) \cap L^2(0,T;H^1(\Omega)), \quad \del_t\rho_i \in L^2(0,T;W^{-1,2}(\Omega))
    \]
    \[
    \theta \in L^\infty(0,T;L^{\alpha+2}(\Omega)), \quad \log\theta \in L^\infty(0,T;L^1(\Omega))
    \]
    \[
    \theta^\frac{\alpha}{2} \in L^2(0,T;H^1(\Omega)), \quad \theta^{\alpha+1} \in L^2(0,T;H^1(\Omega)), \quad \del_t\theta \in L^2(0,T;W^{-1,2}(\Omega)).
    \]
    Furthermore, the solution satisfies the entropy inequality
    \[
    H(\tau)+\int_0^\tau\int_\Omega \mathcal D_{\rm diff}(\rho,J)\,dxdt+\frac4{\alpha^2} \int_0^\tau\int_\Omega |\nabla\theta^{\alpha/2}|^2\,dxdt \le H(0)
    \]
    for almost every $\tau\in(0,T)$, where
    \[
    \mathcal D_{\rm diff}(\rho,J)=\frac12\sum_{i=1}^n\sum_{j\ne i}\frac1{D_{ij}}\left|\sqrt{\frac{\rho_j}{\rho_i}}J_i-\sqrt{\frac{\rho_i}{\rho_j}}J_j\right|^2.
    \]
\end{theorem}

\begin{theorem}[Entropy production defect measure]\label{defmeas}
Let $(\rho_i^{m},\theta^{m},J_i^m)=(\rho_i^{\sigma_m,\epsilon_m},\theta^{\sigma_m,\epsilon_m},J_i^{\sigma_m,\epsilon_m})$ be the approximate solutions used in the proof of the existence theorem, where $\sigma_m,\epsilon_m\to0$. Define the approximate entropy-production measures
\[
d\mu_m = \left[ \mathcal D_{\rm diff}(\rho^m,J^m)+\frac4{\alpha^2}|\nabla(\theta^m)^\frac\alpha2|^2+\epsilon_m|\nabla\log\theta^m|^2+\epsilon_m|\nabla\theta^m|^2\right]dxdt .
\]
Then $(\mu_m)_m$ is bounded in $\mathcal M^+(\overline\Omega\times[0,T])$. Hence, after extracting a subsequence,
\[
\mu_m \stackrel{*}{\rightharpoonup}\mu \quad\textnormal{in } \mathcal M^+(\overline\Omega\times[0,T]).
\]
Moreover,
\[
\mu \ge \left[ \mathcal D_{\rm diff}(\rho,J)+\frac4{\alpha^2}|\nabla\theta^\frac\alpha2|^2 \right]dxdt.
\]
Consequently there exists a nonnegative Radon measure $\mathfrak m\in \mathcal M^+(\overline\Omega\times[0,T])$ such that
\[
\mu = \left[ \mathcal D_{\rm diff}(\rho,J) + \frac4{\alpha^2}|\nabla\theta^\frac\alpha2|^2 \right]dxdt + \mathfrak m .
\]
The measure $\mathfrak m$ is the entropy-production defect measure.

For every $\tau\in(0,T)$ which is a continuity point of the measure $\mu$, and for which the entropy is lower semicontinuity, one has
\[
H(\tau)+\int_0^\tau\int_\Omega\mathcal D_{\rm diff}(\rho,J)\,dxdt+\frac4{\alpha^2} \int_0^\tau\int_\Omega|\nabla\theta^\frac\alpha2|^2\,dxdt+\mathfrak m(\overline\Omega\times(0,\tau))\leq H(0).
\]
In particular, this holds for almost every $\tau\in(0,T)$.
\end{theorem}

\begin{remark}
    In contrast to the isothermal Maxwell--Stefan system, where weak solutions were shown to automatically satisfy the entropy identity \cite{BGT25}, the degenerate nonisothermal system may lose entropy dissipation in the limit.
\end{remark}


\section{Proof of Theorem \ref{thmex}}

Our proof builds upon the existence theory of \cite[Theorem 1]{GJ24}. We therefore introduce an approximate problem that satisfies the hypotheses of \cite[Theorem 1]{GJ24}, which yields the existence of approximate solutions. We then derive estimates that are uniform with respect to the approximation parameters and use them to pass to the limit.

To be compatible with the framework of \cite[Theorem 1]{GJ24}, the approximate problem must satisfy the following assumptions:

\begin{itemize}
    \item[(1)] $\Omega\subset\mathbb{R}^3$ a bounded domain with Lipschitz boundary,
    \item[(2)] $\rho_i^0 \in L^\infty(\Omega)$ with $\rho_i^0\geq0$ for all $i\in\{1,\dots,n\}$ and $0<\rho_*\leq \sum_{i=1}^n\rho_i^0 \leq \rho^*$ for some $\rho_*,\rho^*>0$,
    \item[(3)] $\theta^0 \in L^\infty(\Omega)$, with $\inf_\Omega \theta^0(x)>0$,
    \item[(4)] $c(1+\theta^2) \leq \kappa(\theta) \leq C(1+\theta^2)$, for some $c,C>0$. 
\end{itemize}
\subsection*{The approximate problem} Consider \eqref{weakmass} and \eqref{weakenergy} with initial data

\begin{equation}
    \rho_i^\sigma(x,0) = \frac{\rho_i^0(x)+\sigma}{1+n\sigma}, \quad \theta^\sigma(x,0) = \min\left\{\max\left\{\theta^0(x),\sigma\right\},\frac1\sigma\right\}
\end{equation}
and the regularized thermal conductivity
\begin{equation}
    \kappa^\epsilon(\theta) = \theta^\alpha+\epsilon(1+\theta^2)
\end{equation}
for some $\sigma,\epsilon>0$. Then, assumptions (1) through (4) are satisfied and thus there exists a weak solution $\left(\left(\rho_i^{\sigma,\epsilon}\right)_{i=1}^n, \theta^{\sigma,\epsilon}\right)$, such that
\begin{align}
    & \sqrt{\rho_i^{\sigma,\epsilon}} \in L^\infty(\Omega\times(0,T)) \cap L^2(0,T;H^1(\Omega)), \quad \del_t \rho_i^{\sigma,\epsilon} \in L^2(0,T;W^{-1,2}(\Omega)) \label{estimatedensity} \\
    & \theta^{\sigma,\epsilon} \in L^\infty(0,T;L^2(\Omega)) \cap L^2(0,T;H^1(\Omega)), \quad \del_t\theta^{\sigma,\epsilon} \in L^\frac{16}{11}(0,T;W^{-1,\frac{16}{11}}(\Omega)). \label{estimatetemper}
\end{align}
Estimates \eqref{estimatedensity}--\eqref{estimatetemper} are in principle not uniform in $\sigma$ and $\epsilon$ and need to be treated with caution.

Moreover, the weak solution satisfies the entropy inequality
\begin{equation}\label{entreq}
     \frac{d}{dt} H^{\sigma,\epsilon}(t) + \int_\Omega\mathcal D_{\rm diff}^{\sigma,\epsilon} \,dx + \int_\Omega\mathcal D_{\rm heat}^{\sigma,\epsilon} \,dx \leq 0,
\end{equation}
where
\[
    \mathcal D_{\rm diff}^{\sigma,\epsilon} = \sum_{i=1}^n\sum_{j\not=i} \frac{1}{2D_{ij}}\left|\sqrt{\frac{\rho_j^{\sigma,\epsilon}}{\rho_i^{\sigma,\epsilon}}} J_i^{\sigma,\epsilon} - \sqrt{\frac{\rho_i^{\sigma,\epsilon}}{\rho_j^{\sigma,\epsilon}}} J_j^{\sigma,\epsilon} \right|^2
\]
and
\[
    \mathcal D_{\rm heat}^{\sigma,\epsilon} =  \kappa^\epsilon(\theta^{\sigma,\epsilon})\frac{|\nabla\theta^{\sigma,\epsilon}|^2}{(\theta^{\sigma,\epsilon})^2}.
\]
We note that as $\rho_i\to0$, $J_i=\mathcal{O}(\sqrt{\rho_i})$ and for this reason the dissipation remains well-defined even near vacuum.

Assumptions (3) and (4) are imposed only at the level of the approximate problem. In particular, the boundedness above and below of $\theta^0$ and the quadratic growth condition on $\kappa$ will be removed in the limit.

We note that we restrict our attention to exponents $\alpha\in(0,2]$ so that assumption (4) is satisfied. Indeed, for every $\alpha>0$, $\kappa^\epsilon(\theta)\ge c(1+\theta^2)$, with $c>0$ depending on $\epsilon$, whereas the upper bound $\kappa^\epsilon(\theta)\le C(1+\theta^2)$ holds only when $\alpha\le2$. Although the admissible range of exponents could potentially be extended beyond $2$, doing so would require developing a new existence theory. Since our primary interest is the physically relevant case $\alpha=\frac12$, the restriction $\alpha\in(0,2]$ is sufficient for the purposes of this paper.

Integrating the entropy inequality \eqref{entreq} one obtains two types of estimates:

\noindent (1) the ones coming from $\mathcal D_{\rm diff}^{\sigma,\epsilon}$ and are the same as in \cite{HJ21, GJ24}, namely
\[ 
J_i^{\sigma,\epsilon} \in L^2(\Omega\times(0,T))
\]
which in turn implies 
\[ 
\nabla\sqrt{\rho_i^{\sigma,\epsilon}} \in L^2(\Omega\times(0,T))
\]
and 
\[
\del_t\rho_i^{\sigma,\epsilon} \in L^2(0,T;W^{-1,2}(\Omega))
\]
as well as the uniform bounds $0 \leq \rho_i^{\sigma,\epsilon} \leq 1$. Due to our choice of initial data
\[ 
\int_\Omega |\log(\theta^0)^\sigma| \leq \int_\Omega |\log \theta^0|<+\infty
\]
and so $H^\sigma(0) <+\infty$ independently of $\sigma,\epsilon$, hence the above estimates are now uniform in $\sigma$ and $\epsilon$. 

Combining these estimates, we find
\[
\rho_i^{\sigma,\epsilon} \in L^2(0,T;H^1(\Omega))
\] 
and so, by the Aubin--Lions lemma, $\rho_i^{\sigma,\epsilon} \to \rho_i$ strongly in $L^2(\Omega\times(0,T))$ and up to a subsequence $\rho_i^{\sigma,\epsilon} \to \rho_i$ almost everywhere in $\Omega\times(0,T)$. Since the regularized initial data satisfy $(\rho^0_i)^\sigma \to \rho_i^0$ strongly in $L^2(\Omega)$, we may pass to the limits $\sigma,\epsilon\to0$ in the mass equations.

\noindent (2) the ones coming from $\mathcal D_{\rm heat}^{\sigma,\epsilon}$, namely

\begin{align*}
    \int_0^T \int_\Omega \kappa^\epsilon(\theta^{\sigma,\epsilon}) \frac{|\nabla\theta^{\sigma,\epsilon}|^2}{(\theta^{\sigma,\epsilon})^2} \, dx dt & = \frac{4}{\alpha^2}\int_0^T \int_\Omega  |\nabla(\theta^{\sigma,\epsilon})^\frac{\alpha}{2}|^2 \, dx dt  \\
    & \quad + \epsilon \int_\Omega |\nabla \log\theta^{\sigma,\epsilon}|^2 \, dx dt + \epsilon \int_\Omega |\nabla \theta^{\sigma,\epsilon}|^2 \, dx dt
\end{align*}
and so we obtain the uniform-in-$(\sigma,\epsilon)$ estimate
\begin{equation}\label{entrest}
    \nabla(\theta^{\sigma,\epsilon})^\frac\alpha2 \in L^2(\Omega\times(0,T))
\end{equation}
and the $(\sigma,\epsilon)$-dependent estimates
\begin{equation} \label{nonunifest}
    \sqrt{\epsilon} \nabla \log\theta^{\sigma,\epsilon}, \sqrt{\epsilon} \nabla \theta^{\sigma,\epsilon} \in L^2(\Omega\times(0,T)).
\end{equation}
Notice that in our case the logarithmic estimate does not survive the limit $\epsilon\to0$. Instead, we obtain a much weaker estimate, namely the one in \eqref{entrest}.  

\subsection*{Energy estimates} However, none of these is sufficient to define the heat flux term $J_e = -\kappa(\theta) \nabla\theta$ in the weak formulation of the energy equation. For this reason, better estimates are needed. 

Testing formally against $(\theta^{\sigma,\epsilon})^{\alpha+1}$ we find
\begin{equation} \label{renormenergy}
    \begin{split}
    \int_\Omega (\theta^{\sigma,\epsilon})^{\alpha+2} \,dx & + \frac{\alpha+2}{\alpha+1} \int_0^t \int_\Omega |\nabla (\theta^{\sigma,\epsilon})^{\alpha+1}|^2 \,dxdt \\
    & + \frac{4\epsilon(\alpha+1)}{\alpha+2} \int_0^t \int_\Omega |\nabla (\theta^{\sigma,\epsilon})^{\frac\alpha2+1}|^2 \,dxdt \\
    & + \frac{4\epsilon(\alpha+1)(\alpha+2)}{(\alpha+4)^2} \int_0^t \int_\Omega |\nabla (\theta^{\sigma,\epsilon})^{\frac\alpha2+2}|^2 \,dxdt = \int_\Omega ((\theta^\sigma)^0)^{\alpha+2} \,dx
    \end{split}
\end{equation}
The estimates coming from this renormalized energy equation would be enough to define the heat flux and pass to the limit in the weak formulation of the energy equation, however the quantity $\theta^{\alpha+1}$ is not a valid test function. 

Yet, we may justify the preceding computation by a renormalized testing argument. Fix $k>0$, and set
\[
\beta_k(s):=T_k(s^{\alpha+1}), \qquad B_k(s):=\int_0^s \beta_k(r)\,dr
\]
where $T_k(s)=\min\{s, k\}$. Since $\beta_k$ is bounded, Lipschitz and nondecreasing, and since
\[
0\le B_k(s)\le ks \qquad \text{for } s\ge 0,
\]
the bound $\theta^{\sigma,\epsilon}\in L^\infty(0,T;L^2(\Omega))$ coming from \eqref{estimatetemper} implies, because $\Omega$ is bounded, that
\[
B_k(\theta^{\sigma,\epsilon})\in L^\infty(0,T;L^1(\Omega)).
\]
Moreover, for fixed $\epsilon>0$, the estimate $\sqrt{\epsilon}\,\nabla\theta^{\sigma,\epsilon}\in L^2(\Omega\times(0,T))$ from \eqref{nonunifest} gives
\[
\nabla\theta^{\sigma,\epsilon}\in L^2(\Omega\times(0,T)).
\]
Hence, by the Sobolev chain rule,
\[
\nabla \beta_k(\theta^{\sigma,\epsilon}) = \beta_k'(\theta^{\sigma,\epsilon})\nabla\theta^{\sigma,\epsilon} \qquad \text{a.e. in } \Omega\times(0,T),
\]
and therefore
\[
\beta_k(\theta^{\sigma,\epsilon})\in L^2(0,T;H^1(\Omega)).
\]
Furthermore,
\[
\beta_k'(s) = (\alpha+1)s^\alpha \mathbf 1_{\{s^{\alpha+1}<k\}} \qquad \text{for a.e. } s\ge0.
\]
Thus $\beta_k'$ is supported in the set where
\[
(\theta^{\sigma,\epsilon})^{\alpha+1}<k,
\]
and consequently $\theta^{\sigma,\epsilon}$ is bounded on the support of $\beta_k'(\theta^{\sigma,\epsilon})$. It follows that
\[
\kappa^\epsilon(\theta^{\sigma,\epsilon}) \beta_k'(\theta^{\sigma,\epsilon}) |\nabla\theta^{\sigma,\epsilon}|^2 \in L^1(\Omega\times(0,T)).
\]
Therefore $\beta_k(\theta^{\sigma,\epsilon})$ is an admissible renormalized test function. Applying the weak integration by parts formula obtained by Steklov time regularization for parabolic Neumann problems, as in \cite[Lemmas~3.2--3.4]{A26}, with
\[
w=\theta^{\sigma,\epsilon}, \qquad G=\beta_k, \qquad \widetilde G=B_k,
\]
and then approximating the time cut-off $\mathbf 1_{[0,t]}$, we obtain, for almost every $t\in(0,T)$,
\[
\int_\Omega B_k(\theta^{\sigma,\epsilon}(t))\,dx + \int_0^t\int_\Omega \kappa^\epsilon(\theta^{\sigma,\epsilon}) \beta_k'(\theta^{\sigma,\epsilon}) |\nabla\theta^{\sigma,\epsilon}|^2\,dxds
\le
\int_\Omega B_k((\theta^0)^\sigma)\,dx .
\]
The homogeneous Neumann condition is encoded in the weak formulation, and therefore no boundary contribution appears.

Strictly speaking, if the definition of renormalized solution is formulated only for renormalizations whose derivative has compact support, one introduces an additional cut-off. Namely, choose $\chi_\ell\in W^{1,\infty}([0,\infty))$ such that $0\le \chi_\ell\le1$, $\chi_\ell\equiv1$ on $[0,\ell]$, and $\chi_\ell\equiv0$ on $[2\ell,\infty)$, and define
\[
B_{k,\ell}(s):=\int_0^s \beta_k(v)\chi_\ell(v)\,dv .
\]
Then $B_{k,\ell}'=\beta_k\chi_\ell$ has compact support, so the renormalized formulation applies to $B_{k,\ell}$. Letting $\ell\to\infty$ gives the desired inequality for $B_k$. The passage to the limit follows from monotone convergence for the $B_{k,\ell}$ term and from dominated convergence for the dissipation term, since $\beta_k'$ is supported where $\theta^{\sigma,\epsilon}$ is bounded.

Since
\[
\kappa^\epsilon(s)=s^\alpha+\epsilon(1+s^2),
\]
we have
\[
\kappa^\epsilon(s)\beta_k'(s) = (\alpha+1) \left(s^{2\alpha}+\epsilon s^\alpha +\epsilon s^{\alpha+2}\right) \mathbf 1_{\{s^{\alpha+1}<k\}}
\]
and thus we obtain
\begin{equation} \label{prelim}
\begin{split}
    \int_\Omega B_k(\theta^{\sigma,\epsilon}(t)) \,dx & + \frac{1}{\alpha+1} \int_0^t \int_{\{(\theta^{\sigma,\epsilon})^{\alpha+1} < k\}} |\nabla(\theta^{\sigma,\epsilon})^{\alpha+1}|^2 \,dxdt \\
    & + \frac{4\epsilon(\alpha+1)}{(\alpha+2)^2} \int_0^t \int_{\{(\theta^{\sigma,\epsilon})^{\alpha+1} < k\}} |\nabla (\theta^{\sigma,\epsilon})^{\frac\alpha2+1}|^2 \,dxdt \\
    & + \frac{4\epsilon(\alpha+1)}{(\alpha+4)^2} \int_0^t \int_{\{(\theta^{\sigma,\epsilon})^{\alpha+1} < k\}} |\nabla(\theta^{\sigma,\epsilon})^{\frac\alpha2+2}|^2 \,dxdt \leq \int_\Omega B_k((\theta^{\sigma})^0) \,dx.
\end{split}
\end{equation}
We note that $B_k$ satisfies
\begin{equation} \label{dom}
    B_k \leq B_{k+1} \quad \textnormal{ and } \quad 0 \leq B_k(s) \leq \frac{s^{\alpha+2}}{\alpha+2}
\end{equation}
and since $\theta^0 \in L^{\alpha+2}(\Omega)$, the right-hand side of \eqref{prelim} is uniformly bounded and well-defined
\[
\int_\Omega B_k((\theta^{\sigma})^0) \,dx \leq \frac{1}{\alpha+2} \int_\Omega ((\theta^{\sigma})^0)^{\alpha+2} \,dx.
\]

We, now, pass to the limit $k\to\infty$ on the left-hand side of \eqref{prelim} using the Monotone Convergence Theorem. Indeed, \eqref{dom} implies
\[
B_k(\theta^{\sigma,\epsilon}) \nearrow \frac{(\theta^{\sigma,\epsilon})^{\alpha+2}}{\alpha+2} \quad \textnormal{ as } k\to\infty,
\]
and similarly setting $f_k=|\nabla \left((\theta^{\sigma,\epsilon})^{\alpha+1}\right)|^2 \mathbf 1_{(\theta^{\sigma,\epsilon})^{\alpha+1} < k}$, we have
\[
f_k\geq0, \quad f_k \leq f_{k+1} \quad \textnormal{ and } \quad f_k \nearrow |\nabla (\theta^{\sigma,\epsilon})^{\alpha+1}|^2 \quad \textnormal{ as } k\to\infty.
\]
The other two integrals are treated similarly.

Putting everything together, we obtain \eqref{renormenergy} as an inequality and thus the uniform-in-$(\sigma,\epsilon)$ estimates
\begin{equation}\label{energyest}
    \theta^{\sigma,\epsilon} \in L^\infty(0,T;L^{\alpha+2}(\Omega)), \quad \nabla (\theta^{\sigma,\epsilon})^{\alpha+1} \in L^2(\Omega\times(0,T))
\end{equation}
and the $(\sigma,\epsilon)$-dependent estimates
\begin{equation}\label{nonuniest}
    \sqrt{\epsilon}\nabla (\theta^{\sigma,\epsilon})^{\frac\alpha2+1}, \sqrt{\epsilon}\nabla (\theta^{\sigma,\epsilon})^{\frac\alpha2+2} \in L^2(\Omega\times(0,T)).
\end{equation}

\subsection*{Compactness} Since $\theta^{\sigma,\epsilon} \in L^\infty(0,T;L^{\alpha+2}(\Omega))$, we also have $\theta^{\sigma,\epsilon} \in L^\infty(0,T;L^{\frac\alpha2}(\Omega))\cap L^\infty(0,T;L^{\alpha+1}(\Omega))$. In particular, we have control over the means of $(\theta^{\sigma,\epsilon})^\frac\alpha2$ and $(\theta^{\sigma,\epsilon})^{\alpha+\frac12}$ in $\Omega$, which combined with the gradient estimates, give by the Poincare--Wirtinger inequality, $(\theta^{\sigma,\epsilon})^{\alpha+1}, (\theta^{\sigma,\epsilon})^{\frac\alpha2} \in L^2(\Omega\times(0,T))$, and so we obtain that $(\theta^{\sigma,\epsilon})^{\frac\alpha2}$ and $(\theta^{\sigma,\epsilon})^{\alpha+1}$ are bounded in $L^2(0,T;H^1(\Omega))$.

Let us, now, identify the space where $J_e^{\sigma,\epsilon}$ lies. We write
\[
J_e^{\sigma,\epsilon} = -\kappa^\epsilon(\theta^{\sigma,\epsilon}) \nabla\theta^{\sigma,\epsilon} = -\frac{1}{\alpha+1} \nabla (\theta^{\sigma,\epsilon})^{\alpha+1} - \epsilon\nabla\theta^{\sigma,\epsilon} - \frac{\epsilon}{3}\nabla(\theta^{\sigma,\epsilon})^3.
\]
The first two terms are in $L^2(0,T;L^2(\Omega))$ because of \eqref{energyest} and \eqref{nonunifest} respectively. The last one is more delicate. In particular, we have
\[
\epsilon \theta^2 \nabla\theta = \frac{2\sqrt{\epsilon}}{\alpha+4} \theta^{1-\frac\alpha2} (\sqrt{\epsilon}\nabla \theta^{\frac\alpha2+2})
\]
where the exponent of $\theta$ is nonnegative because of $\alpha\leq2$ and $\sqrt{\epsilon} \nabla \theta^{\frac\alpha2+2}$ is in $L^2(\Omega\times(0,T))$. 

It is interesting that this particular estimate needs the assumption $\alpha \leq2$, which is the same assumption we needed to make sure that the thermal conductivity $\kappa^\epsilon$ satisfies the bounds $c(1+\theta^2) \leq \kappa^\epsilon(\theta) \leq C(1+\theta^2)$, so that we may apply the existence theory of \cite{GJ24}.

Now, since $\theta^{\sigma,\epsilon} \in L^\infty(0,T;L^{\alpha+2}(\Omega))$, it follows that $(\theta^{\sigma,\epsilon})^{1-\frac\alpha2}\in L^\infty(0,T;L^\frac{\alpha+2}{1-\frac\alpha2}(\Omega))$ and thus $\epsilon \theta^2 \nabla\theta \in L^r(\Omega\times(0,T))$, with
\[
\frac{1}{r} = \frac{1}{2} + \frac{1-\frac\alpha2}{\alpha+2}, \quad \textnormal{ i.e. } \quad r=\frac{\alpha}{2}+1.
\]
Therefore, $\epsilon \theta^2 \nabla\theta \to 0$ strongly in $L^r(\Omega\times(0,T))$ and due to $\alpha\leq2$, it follows $r\leq2$ and thus $J_e^{\sigma,\epsilon} \in L^r(\Omega\times(0,T))$. Hence,
\[
\del_t\theta^{\sigma,\epsilon} = -\diver J_e^{\sigma,\epsilon} \in L^r(0,T;W^{-1,r}(\Omega)), \quad r=\frac{\alpha+2}{2}>1.
\]

Combining $(\theta^{\sigma,\epsilon})^{\alpha+1} \in L^2(0,T;H^1(\Omega))$ with $\del_t\theta^{\sigma,\epsilon} \in L^r(0,T;W^{-1,r}(\Omega))$ we may apply the Aubin--Lions--Dubinskii compactness lemma in the version of \cite{CJL14} to derive strong compactness of $(\theta^{\sigma,\epsilon})_{\sigma,\epsilon}$. Indeed, using the notation of \cite{CJL14}, we may set $m=\alpha+1>1$, $p=m+1$ and the spaces $B=L^p(\Omega)$, $Y=W^{-1,r}(\Omega)$ and $M_+=\{w\geq0: w^m \in H^1(\Omega)\}$, with $[w]_{M_+}=\|w^m\|_{H^1(\Omega)}^{\frac1m}$. The compact embedding of $M_+ \subset\subset B$ follows because if $w_\ell^m$ is bounded in $H^1(\Omega)$, then in space dimension 3, $w_\ell^m$ is compact in $L^q(\Omega)$ for every $q\in[1,6)$. Choose $q=\frac{m+1}{m}=\frac{\alpha+2}{\alpha+1}<6$. Then $w_\ell^m \to z$ strongly in $L^\frac{m+1}{m}(\Omega)$ and hence $w_\ell = (w_\ell^m)^\frac{1}{m}\to z^\frac1m$ strongly in $L^p(\Omega)$. Moreover, $B=L^p(\Omega) \hookrightarrow W^{-1,r}(\Omega)=Y$. Now, our estimates imply boundedness in $L^{2m}(0,T;M_+)$. Indeed,
\[
\|\theta^{\sigma,\epsilon}\|_{L^{2m}(0,T;M_+)}^{2m} = \int_0^T \|(\theta^{\sigma,\epsilon})^m(t)\|^2_{H^1(\Omega)} \,dt
\]
and we already have $(\theta^{\sigma,\epsilon})^{\alpha+1}$ is bounded in $L^2(0,T;H^1(\Omega))$. Since also $\del_t\theta^{\sigma,\epsilon} \in L^r(0,T;Y)$, the Aubin--Lions--Dubinskii compactness lemma yields $\theta^{\sigma,\epsilon} \to \theta$ strongly in $L^{2m}(0,T;L^p(\Omega))=L^{2(\alpha+1)}(0,T;L^{\alpha+2}(\Omega))$.

It remains to identify the nonlinear weak limit. Since $(\theta^{\sigma,\epsilon})^{\alpha+1}$ is bounded in $L^2(0,T;H^1(\Omega))$, there exists $\chi \in L^2(0,T;H^1(\Omega))$, such that up to a subsequence $(\theta^{\sigma,\epsilon})^{\alpha+1} \rightharpoonup \chi$ weakly in $L^2(0,T;H^1(\Omega))$. Due to the strong convergence of $\theta^{\sigma,\epsilon}$, we also have that $(\theta^{\sigma,\epsilon})^{\alpha+1} \to \theta^{\alpha+1}$ strongly $L^2(0,T;L^\frac{p}{m}(\Omega))$, since $1<\frac{p}{m}<6$. Since also $(\theta^{\sigma,\epsilon})^{\alpha+1} \rightharpoonup \chi$ weakly in $L^2(0,T;H^1(\Omega))$ the two limits must coincide because $H^1(\Omega) \hookrightarrow L^\frac{p}{m}(\Omega)$, i.e. $\chi=\theta^{\alpha+1}$.

Following the same ideas, we identify the weak limit of the time derivative in the duality pairing term. Due to $\del_t \theta^{\sigma,\epsilon} \in L^r(0,T;W^{-1,r}(\Omega))$, there exists $\xi \in L^r(0,T;W^{-1,r}(\Omega))$ such that, up to a subsequence, we have $\del_t \theta^{\sigma,\epsilon} \rightharpoonup \xi$ weakly in $L^r(0,T;W^{-1,r}(\Omega))$. Because of the strong compactness of $\theta^{\sigma,\epsilon}$ in $L^{2(\alpha+1)}(0,T;L^{\alpha+2}(\Omega))$, we have $\theta^{\sigma,\epsilon} \to \theta$ strongly in $L^1(0,T;W^{-1,r}(\Omega))$, since $L^{\alpha+2}(\Omega) \hookrightarrow W^{-1,r}(\Omega)$ in 3 space dimensions. Therefore $\xi=\del_t\theta$ in $\mathcal{D}'(0,T;W^{-1,r}(\Omega))$ and so $\del_t\theta^{\sigma,\epsilon} \rightharpoonup \del_t\theta$ weakly in $L^r(0,T;W^{-1,r}(\Omega))$.

Putting everything together, we may pass to the limits $\sigma,\epsilon\to0$ in the energy equation. In fact, after passing to the limits, we end up with 
\[
J_e = -\frac{1}{\alpha+1} \nabla \theta^{\alpha+1} \in L^2(\Omega\times(0,T))
\]
and so 
\[
\del_t\theta \in L^2(0,T;W^{-1,2}(\Omega)).
\]

\subsection*{Linear system and constraint} It remains to pass to the limit in the linear system and the constraint. For the linear system, the convergence $\rho_i^{\sigma,\epsilon} \to \rho_i$ strongly in $L^2(\Omega\times(0,T))$ and $J_i^{\sigma,\epsilon} \rightharpoonup J_i$ weakly in $L^2(\Omega\times(0,T))$ is enough to pass to the limit in the left-hand side weakly in $L^1$ and thus in distributions. For the right-hand side, the strong compactness of $\rho_i^{\sigma,\epsilon}$ and $\theta^{\sigma,\epsilon}$ implies that $\rho_i^{\sigma,\epsilon}\theta^{\sigma,\epsilon} \to \rho_i\theta$ in $L^1(\Omega\times(0,T))$, and so we may write
\[
\int_0^T \int_\Omega \nabla(\rho_i^{\sigma,\epsilon} \theta^{\sigma,\epsilon}) \cdot \varphi \,dxdt = - \int_0^T \int_\Omega \rho_i^{\sigma,\epsilon} \theta^{\sigma,\epsilon} \diver \varphi \,dxdt \to - \int_0^T \int_\Omega \rho_i \theta \diver \varphi \,dxdt
\]
that is $\nabla(\rho_i^{\sigma,\epsilon} \theta^{\sigma,\epsilon}) \to \nabla(\rho_i \theta)$ in $\mathcal{D}'(\Omega\times(0,T))$. So the linear system is satisfied in the sense of distributions. 

For the constraint, since $\sum_{i=1}^n J_i^{\sigma,\epsilon} = 0$, we have that for all $\varphi \in C^\infty_0(\Omega)$
\[
\int_\Omega \sum_{i=1}^n J_i^{\sigma,\epsilon} \varphi = 0.
\]
Then, the convergence $J_i^{\sigma,\epsilon} \rightharpoonup J_i$ weakly in $L^2(\Omega\times(0,T))$, implies that
\[
\int_\Omega \sum_{i=1}^n J_i \varphi = 0, \quad \forall ~ \varphi \in C^\infty_0(\Omega) \textnormal{ with } \sum_{i=1}^n J_i \in L^2(\Omega\times(0,T)).
\]
Therefore, $\sum_{i=1}^n J_i=0$ almost everywhere in $\Omega\times(0,T)$.

\subsection*{Passage to the limit in the entropy inequality} We now pass to the limit in the entropy inequality. The approximate solutions satisfy, for almost every $\tau\in(0,T)$,
\[
\begin{aligned}
    H^{\sigma,\epsilon}(\tau) &+ \int_0^\tau\int_\Omega \mathcal D_{\rm diff}(\rho^{\sigma,\epsilon},J^{\sigma,\epsilon})\,dxdt + \frac4{\alpha^2} \int_0^\tau\int_\Omega |\nabla(\theta^{\sigma,\epsilon})^\frac{\alpha}2|^2\,dxdt  \\
    &+ \epsilon \int_0^\tau\int_\Omega |\nabla\log\theta^{\sigma,\epsilon}|^2\,dxdt + \epsilon \int_0^\tau\int_\Omega |\nabla\theta^{\sigma,\epsilon}|^2\,dxdt \le H^\sigma(0).
\end{aligned}
\]
The last two terms are nonnegative. Dropping them gives
\[
H^{\sigma,\epsilon}(\tau)+\int_0^\tau\int_\Omega \mathcal D_{\rm diff}(\rho^{\sigma,\epsilon},J^{\sigma,\epsilon})\,dxdt + \frac4{\alpha^2} \int_0^\tau\int_\Omega \nabla(\theta^{\sigma,\epsilon})^\frac{\alpha}2|^2\,dxdt \le H^\sigma(0).
\]
Since $0\le\rho_i^{\sigma,\epsilon}\le1$, the density entropy is uniformly bounded:
\[
-n|\Omega| \le \int_\Omega \sum_{i=1}^n \rho_i^{\sigma,\epsilon} (\log\rho_i^{\sigma,\epsilon}-1)\,dx \le0 .
\]
Moreover,
\[
\log s \le C(1+s^{\alpha+2}) \qquad\text{for all }s\ge0.
\]
Therefore, the strong convergence of $\rho_i^{\sigma,\epsilon}$ and $\theta^{\sigma,\epsilon}$ implies
\[
H(\tau) \le \liminf_{\sigma,\epsilon\to0} H^{\sigma,\epsilon}(\tau).
\]

The diffusive part of the entropy production is lower semicontinuous with respect to
\[
\rho_i^{\sigma,\epsilon}\to\rho_i \quad\text{strongly}, \qquad J_i^{\sigma,\epsilon}\rightharpoonup J_i \quad\text{weakly}.
\]
Therefore
\[
\int_0^\tau\int_\Omega \mathcal D_{\rm diff}(\rho,J)\,dxdt \le \liminf_{\sigma,\epsilon\to0} \int_0^\tau\int_\Omega \mathcal D_{\rm diff}(\rho^{\sigma,\epsilon},J^{\sigma,\epsilon}) \,dxdt.
\]
Similarly, since
\[
(\theta^{\sigma,\epsilon})^\frac\alpha2 \rightharpoonup \theta^\frac\alpha2 \quad\text{weakly in }L^2(0,T;H^1(\Omega)),
\]
we have
\[
\int_0^\tau\int_\Omega |\nabla\theta^\frac\alpha2|^2\,dxdt \le \liminf_{\sigma,\epsilon\to0} \int_0^\tau\int_\Omega |\nabla(\theta^{\sigma,\epsilon})^\frac\alpha2|^2\,dxdt.
\]

For the initial entropy, since
\[
\log\theta_0^\sigma = \min\{\max\{\log\theta^0,\log\sigma\},-\log\sigma\},
\]
we have
\[
\log\theta_0^\sigma\to\log\theta^0 \quad\text{strongly in }L^1(\Omega).
\]
Indeed, this is just the convergence of symmetric truncations of the $L^1$-function $\log\theta^0$. Therefore
\[
\int_\Omega -\log\theta_0^\sigma\,dx \to \int_\Omega -\log\theta^0\,dx .
\]
The density part of the initial entropy converges as well, since $\rho_i^{0,\sigma}\to\rho_i^0$ strongly in $L^p(\Omega)$ for every finite $p$, $0\le\rho_i^{0,\sigma}\le1$, and $s\mapsto s(\log s-1)$ is continuous and bounded on $[0,1]$. Consequently,
\[
H^\sigma(0)\to H^0 \quad\text{as }\sigma\to0.
\]
Putting everything together, we obtain 
\[
H(\tau)+\int_0^\tau\int_\Omega \mathcal D_{\rm diff}(\rho,J)\,dxdt+\frac4{\alpha^2} \int_0^\tau\int_\Omega|\nabla\theta^\frac\alpha2|^2\,dxdt \le H(0)
\]
for almost every $\tau\in(0,T)$. Therefore, the limiting solution has finite entropy and satisfies the entropy inequality. In particular, 
\[
\log\theta \in L^\infty(0,T;L^1(\Omega))
\]
that is, $\theta>0$ almost everywhere in $\Omega\times(0,T)$. Indeed, if $\theta=0$ on a set of positive measure, then $\log\theta=-\infty$ on that set, contradicting the above bound.


\section{Proof of Theorem \ref{defmeas}}

Let $(\rho_i^m,\theta^m,J_i^m)=(\rho_i^{\sigma_m,\epsilon_m},\theta^{\sigma_m,\epsilon_m},J_i^{\sigma_m,\epsilon_m})$ be the approximate sequence extracted in the existence proof, with $\sigma_m,\epsilon_m\to0$. Define
\[
d\mu_m=\left[\mathcal D_{\rm diff}(\rho^m,J^m)+\frac4{\alpha^2}|\nabla(\theta^m)^\frac\alpha2|^2+\epsilon_m|\nabla\log\theta^m|^2+\epsilon_m|\nabla\theta^m|^2\right]dxdt .
\]
The approximate entropy inequality and the convergence of the initial entropy imply
\[
\sup_m \mu_m(\overline\Omega\times[0,T])<\infty .
\]
Therefore, by the weak-* compactness of bounded sets in the space of nonnegative Radon measures, there exists a subsequence, not relabeled, and a measure $\mu\in\mathcal M^+(\overline\Omega\times[0,T])$ such that
\[
\mu_m\stackrel{*}{\rightharpoonup}\mu \quad\text{in } \mathcal M^+(\overline\Omega\times[0,T]).
\]

We now identify the absolutely continuous part of $\mu$. Let $\varphi\in C(\overline\Omega\times[0,T])$, $\varphi\ge0$. Since
\[
(\theta^m)^\frac\alpha2 \rightharpoonup \theta^\frac\alpha2 \quad\text{weakly in }L^2(0,T;H^1(\Omega)),
\]
we have
\[
\sqrt{\varphi}\,\nabla(\theta^m)^\frac\alpha2 \rightharpoonup \sqrt{\varphi}\,\nabla\theta^\frac\alpha2 \quad\text{weakly in }L^2(\Omega\times(0,T)).
\]
Hence
\[
\int_0^T\int_\Omega\varphi|\nabla\theta^\frac\alpha2|^2\,dxdt \leq \liminf_{m\to\infty} \int_0^T\int_\Omega\varphi |\nabla(\theta^m)^\frac\alpha2|^2\,dxdt.
\]
Multiplying by $\frac4{\alpha^2}$, we obtain the lower semicontinuity of the heat part.

Similarly, by the lower semicontinuity of the diffusive part under the convergence
\[
\rho^m\to\rho \quad\text{strongly}, \qquad J^m\rightharpoonup J \quad\text{weakly},
\]
one has
\[
\int_0^T\int_\Omega\varphi\, \mathcal D_{\rm diff}(\rho,J)\,dxdt \le \liminf_{m\to\infty} \int_0^T\int_\Omega\varphi\, \mathcal D_{\rm diff} (\rho^m,J^m)\,dxdt .
\]
The remaining two terms,
\[
\epsilon_m|\nabla\log\theta^m|^2, \qquad \epsilon_m|\nabla\theta^m|^2,
\]
are nonnegative and therefore can only add mass to the weak-* limit. Consequently,
\[
\int_0^T\int_\Omega\varphi \left[\mathcal D_{\rm diff}(\rho,J)+\frac4{\alpha^2}|\nabla\theta^\frac\alpha2|^2\right]\,dxdt \le \int_0^T\int_{\overline{\Omega}}\varphi\,d\mu .
\]
Since this holds for every nonnegative $\varphi\in C(\overline{\Omega}\times(0,T))$, we obtain the measure inequality
\[
\mu \ge \left[ \mathcal D_{\rm diff}(\rho,J) + \frac4{\alpha^2}|\nabla\theta^\frac\alpha2|^2 \right]dxdt .
\]
Thus the measure
\[
\mathfrak m := \mu - \left[ \mathcal D_{\rm diff}(\rho,J) + \frac4{\alpha^2}|\nabla\theta^\frac\alpha2|^2 \right]dxdt 
\]
is a nonnegative Radon measure on $\overline{\Omega}\times(0,T)$.

Let $\tau\in(0,T)$ be a continuity point of the measure $\mu$ and a time for which the entropy is lower semicontinuous. Then
\[
\mu_m(\overline\Omega\times(0,\tau)) \to \mu(\overline\Omega\times(0,\tau)).
\]
The approximate entropy inequality gives
\[
H^m(\tau) + \mu_m(\overline\Omega\times(0,\tau)) \leq H^m(0).
\]
Passing to the limit, using
\[
H(\tau) \leq \liminf_{m\to\infty} H^m(\tau), \qquad H^m(0)\to H(0),
\]
we obtain
\[
H(\tau) + \mu(\overline\Omega\times(0,\tau)) \leq H(0)
\]
and using the decomposition
\[
\mu = \left[ \mathcal D_{\rm diff}(\rho,J) + \frac4{\alpha^2}|\nabla\theta^\frac\alpha2|^2 \right]dxdt + \mathfrak m ,
\]
we conclude that
\[
H(\tau)+\int_0^\tau\int_\Omega \mathcal D_{\rm diff}(\rho,J)\,dxdt + \frac4{\alpha^2} \int_0^\tau\int_\Omega |\nabla\theta^\frac\alpha2|^2\,dxdt + \mathfrak m(\overline\Omega\times(0,\tau)) \leq H(0).
\]
Since finite Radon measures have at most countably many atoms in time and the lower semicontinuity holds for almost every time, the above inequality holds for almost every $\tau\in(0,T)$. This proves the result.



\begin{thebibliography}{10}

\bibitem{A26} 
\textnormal{\textsc{M. Aoun}, Existence and uniqueness of renormalized solutions for parabolic Neumann problem with $L^1$ data {\it Nonlinear Analysis} {\bf 266}, (2026) 114037.}

\bibitem{AP15} 
\textnormal{\textsc{S. Axmann, M. Pokorn\'y}, Time-periodic solutions to the full Navier--Stokes--Fourier system with radiation on the boundary, {\it J. Math. Anal. Appl.} {\bf 428} (2015), 414-444.}

\bibitem{BGT25} 
\textnormal{\textsc{L.C. Berselli, S. Georgiadis, A.E. Tzavaras}, Absence of anomalous dissipation for the Maxwell--Stefan system, {\it Nonlinearity} {\bf 38} (2025), 025018.}

\bibitem{BD06} 
\textnormal{\textsc{D. Bresch, B. Desjardins}, Stabilité de solutions faibles globales pour les équations
de Navier--Stokes compressible avec température, {\it C. R. Acad. Sci. Paris, Ser. I} {\bf 343} (2006), 219-224.}

\bibitem{BJPZ22} 
\textnormal{\textsc{M. Buli\v{c}ek, A. J\"ungel, M. Pokorn\'y, N. Zamponi}, Existence analysis of a stationary compressible fluid model for heat-conducting and chemically reacting mixtures, {\it J. Math. Phys.} {\bf 63} (2022), 051501.}

\bibitem{CJL14} 
\textnormal{\textsc{X. Chen, A. J\"ungel, J. Liu}, A Note on Aubin--Lions--Dubinskii Lemmas, {\it Acta Appl. Math.} {\bf 133} (2014), 33-43.}

\bibitem{DGZ17} 
\textnormal{\textsc{R. Duan, A. Guo, C. Zhu}, Global strong solution to compressible Navier--Stokes equations with density dependent viscosity and temperature dependent heat conductivity, {\it J. Diff. Eq.} {\bf 262} (2017), 4314-4335.}

\bibitem{FJSZ21} 
\textnormal{\textsc{G. Favre, A. J\"ungel, C. Schmeiser, N. Zamponi}, Existence analysis of a degenerate diffusion system for heat-conducting gases, {\it Nonlinear Differ. Equ. Appl. } {\bf 28} 41 (2021).}

\bibitem{FN17}
\textnormal{\textsc{E. Feireisl, A. Novotn\'y}, Singular Limits in Thermodynamics of Viscous Fluids, Birkh\"auser Cham (2017).}

\bibitem{FNS16} 
\textnormal{\textsc{E. Feireisl, A. Novotn\'y, Y. Sun}, On the motion of viscous, compressible, and heat-conducting liquids, {\it J. Math. Phys.} {\bf 57} (2016), 083101.}

\bibitem{GJ24} 
\textnormal{\textsc{S. Georgiadis, A. J\"ungel}, Global existence of weak solutions and weak--strong uniqueness for nonisothermal Maxwell--Stefan systems, {\it Nonlinearity} {\bf 37} (2024), 075016.}

\bibitem{HJ21} 
\textnormal{\textsc{C. Helmer, A. J\"ungel}, Analysis of Maxwell--Stefan systems for heat conducting fluid mixtures, {\it Nonlin. Anal.: Real World Appl.} {\bf 59} (2021), 103263.}

\bibitem{JS13} 
\textnormal{\textsc{A. J\"ungel, I.V. Stelzer}, Existence analysis of Maxwell-Stefan systems for multicomponent mixtures, {\it SIAM J. Math. Anal.} {\bf 45} (2013), 2421-2440.}

\bibitem{LPP26} 
\textnormal{\textsc{M. Liu, Y. Peng, Z. Peng}, On outer pressure problem of compressible Navier--Stokes system with degenerate heat-conductivity in unbounded domains, {\it J. Diff. Eq.} {\bf 453} (2026), 113909.}

\bibitem{MPZ15} 
\textnormal{\textsc{P.B. Mucha, M. Pokorn\'y, E. Zatorska}, Heat-Conducting, Compressible Mixtures with Multicomponent Diffusion: Construction of a Weak Solution, {\it SIAM J. Math. Anal.} {\bf 47} (2015), 3747-3797.}

\end{thebibliography}
\end{document}